\documentclass[smallextended]{svjour3}     
\smartqed  
\usepackage{mathptmx}      
\usepackage{latexsym}
\usepackage{amsmath}
\usepackage{amssymb}
\usepackage{mathrsfs}
\usepackage{hyperref}

\begin{document}

\title{Generalization of the Glasser -
Manna - Oloa integral and some new integrals of similar type}

\titlerunning{Generalization of the Glasser -
Manna - Oloa integral}

\author{N.M. Vildanov}

\institute{N.M. Vildanov \at
              I.E.Tamm Department of Theoretical Physics, P.N.Lebedev
Physics Institute - 119991 Moscow, Russia \\
              Tel.:none\\
              Fax: none\\
              \email{niyazvil@list.ru}           
}

\maketitle

\begin{abstract}
 As was shown in the previous works by other authors, Glasser - Manna
- Oloa integral arise in the study of the Laplace transform of the
dilogarithm function and can be evaluated in a closed form. In
this article, we give a one parametric generalization of the
Glasser - Manna - Oloa integral. The method employed in the course
of derivation allows to obtain some new integrals of similar type.
These include also integral representation of the Hurwitz zeta
function and some beautiful formulae involving the logarithm of
the gamma function of the argument $-ix+\ln(2\cos x)$.
\end{abstract}

\section{Motivation and results}
\label{intro}

Recently, there has been some interest in the integrals of the form
\begin{equation}\label{oloaa}
    M(a)=\frac{4}{\pi}\int_{0}^{\frac{\pi}{2}}\frac{x^2}{x^2+\ln^2(2e^{-a}\cos x)}d
    x\,.
\end{equation}
\cite{glasser,oloa,amdeberhan,bailey,dixit}. In \cite{glasser}, it
was shown that $M(a)$ is related to the Laplace transform of the
dilogarithm function $\psi(s+1)$
\begin{equation}\label{laplace}
    L(a)=\int_0^\infty e^{-as}\psi(s+1)~ds
\end{equation}
through
\begin{equation}\label{relation}
    L(a)=M(a)-\frac{\gamma}{a}-\frac{\ln(e^a-1)}{1-e^{-a}}H(\ln
    2-a)\,,
\end{equation}
where $H$ is the unit step function and $\gamma$ is the Euler's
constant. The value $M(0)$ can be evaluated in a closed form and
reads
\begin{equation}\label{oloa}
    M(0)=\frac{4}{\pi}\int_{0}^{\frac{\pi}{2}}\frac{x^2}{x^2+\ln^2(2\cos x)}d x=\frac{1}{2}\left(1-\gamma+\ln
    2\pi\right)\,.
\end{equation}

In this article, we give a proof of the following formula
\begin{equation}\label{generalization}
    {\frac{1}{2i}} \int_{-\frac{\pi}{2}}^{\frac{\pi}{2}}\frac{x\left(1+e^{-2ix}\right)^\beta}{\ln\left(1+e^{-2ix}\right)}d x
    =\frac{\pi}{8}\left(1+\ln
    2\pi-\gamma(2\beta+1)-2\ln\Gamma(\beta+1)\right)\,,
\end{equation}
provided ${\rm Re} ~\beta> -1$. This is a one-parametric
generalization of the integral \eqref{oloa}, which directly follows
from \eqref{generalization} when $\beta=0$ \footnote{It is easy to
check that ${\frac{1}{2i}}
\int_{-\frac{\pi}{2}}^{\frac{\pi}{2}}\frac{x}{\ln\left(1+e^{-2ix}\right)}d
x=\int_{0}^{\frac{\pi}{2}}\frac{x^2}{x^2+\ln^2(2\cos x)}d x$}.

Using the method of derivation of \eqref{generalization} the
following integral representation of the Hurwitz zeta function
\begin{equation}\label{hurwitz}
    I(\alpha,\beta)=\int_{-\frac{\pi}{2}}^{\frac{\pi}{2}}\left(1+e^{-2ix}\right)^\beta \ln^\alpha\left(1+e^{-2ix}\right)d x
    =-\frac{\pi}{\Gamma(-\alpha)}\zeta(\alpha+1,\beta+1)
\end{equation}
can be obtained. \eqref{hurwitz} is valid for $\beta> -1$ and
arbitrary complex $\alpha$; if ${\rm Re}~\beta=-1$, then ${\rm
Re}~\alpha < -1$. For $\alpha=-1$, \eqref{hurwitz} reduces to
\begin{equation}\label{alphaminus1}
    \int_{-\frac{\pi}{2}}^{\frac{\pi}{2}}\frac{\left(1+e^{-2ix}\right)^\beta}{\ln\left(1+e^{-2ix}\right)}d x=\frac{\pi}{2}(1+2\beta)
\end{equation}
and for $\beta =0$ it is equivalent to the integral representation
of the Riemann zeta function
\begin{equation}\label{riemann}
    \int_{-\frac{\pi}{2}}^{\frac{\pi}{2}} \ln^{\alpha-1}{\left({1+e^{-2ix}}\right)}d
    x=-\frac{\pi\zeta(\alpha)}{\Gamma(1-\alpha)}\,.
\end{equation}
Integral in \eqref{riemann} converges for arbitrary complex
$\alpha$.

The following integrals can be obtained from \eqref{alphaminus1} and
\eqref{riemann}:
\begin{equation}\label{int1}
    \int_{0}^{\frac{\pi}{2}}\frac{\ln(2\cos x)}{x^2+\ln^2(2\cos x)}~d
    x=\frac{\pi}{4}\,,
\end{equation}
\begin{equation}\label{int3}
    \int_0^{\frac{\pi}{2}}\ln\left[x^2+\ln^2(2\cos x)\right]dx=0\,,
\end{equation}
\begin{equation}\label{int2}
    \int_{0}^{\frac{\pi}{2}}\frac{1}{x^2+\ln^2(2\cos x)}~d
    x-2\int_{0}^{\frac{\pi}{2}}\frac{x^2}{\left[x^2+\ln^2(2\cos x)\right]^2}~d
    x=\frac{\pi}{24}\,.
\end{equation}
Putting $\beta =0$ in \eqref{alphaminus1}, we obtain \eqref{int1}.
\eqref{int3} is obtained from \eqref{riemann} differentiating by
$\alpha$ at $\alpha=1$ (see also \eqref{intlog}). Assuming
$\alpha=-1$ in \eqref{riemann} and taking into account the value of
the Riemann zeta function $\zeta(-1)=-\frac{1}{12}$ yields
\eqref{int2}.

In the light of formula \eqref{int3} it is natural to study the
integrals
\begin{equation}\label{intlog}
    f(a)=\int_0^{\frac{\pi}{2}}\ln\left[x^2+\ln^2(2e^{-a}\cos
    x)\right]dx=\pi\ln \frac{a}{e^b-1}\,,
\end{equation}
\begin{equation}\label{intlog3}
    g(a)=\int_0^{\frac{\pi}{2}}\ln\left[x^2+\ln^2(2e^{-a}\cos
    x)\right]\cos2x~dx=\frac{\pi}{2}\left(1-\frac{1}{a}-e^b+\frac{1}{e^b-1}\right)\,,
\end{equation}
where $b=\min\{ a,\ln2\}$. Surprisingly they can be evaluated for
any real $a$.

From \eqref{intlog}, the following identities can be deduced
\begin{equation}\label{remarkable}
    \int_0^{\frac{\pi}{2}}\ln\left|\frac{\Gamma[1-\frac{ix+\ln(2\cos x)}{\ln r}]}{\Gamma[1-\frac{ix-\ln(2\cos
    x)}{\ln r}]}\right|
    ~dx=\frac{\pi}{2}\ln\prod_{n=1}^\infty(1-r^n)\,,\qquad
    0<r<\frac{1}{2}
\end{equation}
\begin{equation}\label{remarkable2}
    \int_0^{\frac{\pi}{2}}\ln\left|\Gamma[c+ix-\ln(2\cos x)]\right|
    ~dx=\frac{\pi}{2}\ln\Gamma(c)\,,\quad
    c\geq\ln 2
\end{equation}

In the rest of the article we give the proofs of the formulas listed
above, namely, \eqref{generalization},\eqref{hurwitz},
(\ref{intlog}-\ref{remarkable2}).

\section{Proofs of the new results presented in section \ref{intro}}

\paragraph{ Proof of formula \eqref{generalization}.} Let us make the
substitution $y=\ln\left(1+e^{-2ix}\right)$ in the integral
$$
I(\beta)={\frac{1}{2i}}
\int_{-\frac{\pi}{2}}^{\frac{\pi}{2}}\frac{x\left(1+e^{-2ix}\right)^\beta}{\ln\left(1+e^{-2ix}\right)}d
x\,.
$$
The result is
\begin{equation}\label{}
    I(\beta)=\frac{1}{8i}\int_{-\infty}^{(0+)}\frac{\ln\left(e^{y}-1\right)}{y}\frac{e^{y(\beta+1)}}{e^y-1}d
    y\,.
\end{equation}
This integral should be understood in the following sense: there is
a cut along the line ($-\infty, 0$] in the complex $y$ plane; the
path extends from $-\infty$, circumvents the origin in the
anticlockwise direction and goes back to $-\infty$ on the opposite
side of the cut (Hankel contour; see, e.g., \cite{whittaker}). It is
also assumed that the path does not contain the points $y=\pm 2\pi i
n$, ~$n=1,2,3,...$, which are the poles of $\frac{1}{e^y-1}$. To be
more specific, we will assume that the contour is composed of three
parts: $L_1$ -- the line $(-\infty,\varepsilon)$ in the lower side
of the cut, the circle $\Gamma_\varepsilon$ of radius
$\varepsilon\to 0$ traversed in the anticlockwise direction and
$L_2$ -- the line $(-\infty,\varepsilon)$ on the upper side of the
cut. Such deformation of the path is allowed, since no poles are
intercepted. The logarithm is
$~\ln\left(e^{y}-1\right)=\ln\left|e^{y}-1\right|\mp \pi i$ on the
lower and upper sides of the cut, respectively.

Let us calculate the following two integrals
\begin{equation}\label{}
    I_L(\beta)=\frac{1}{8i}\int_{L_1+L_2}\frac{\ln\left(e^{y}-1\right)}{y}\frac{e^{y(\beta+1)}}{e^y-1}d
    y\,,
\end{equation}
\begin{equation}\label{}
    I_{\Gamma_\varepsilon}(\beta)=\frac{1}{8i}\int_{\Gamma_\varepsilon}\frac{\ln\left(e^{y}-1\right)}{y}\frac{e^{y(\beta+1)}}{e^y-1}d
    y
\end{equation}
separately. First, we rewrite $I_L(\beta)$ as follows
\begin{multline}\label{rewrite}
    I_L(\beta)=-\frac{\pi}{4}\int_{\varepsilon}^{\infty}\frac{e^{-\beta y}
    }{y(e^y-1)}d y\\
    =-\frac{\pi}{4}\int_{\varepsilon}^{\infty}\frac{d
    y}{y}\left\{\frac{e^{-\beta
    y}}{e^y-1}-\frac{1}{y\left[1+y(\frac{1}{2}+\beta)\right]}\right\}
    -\frac{\pi}{4}\int_{\varepsilon}^{\infty}\frac{d
    y}{y^2\left[1+y(\frac{1}{2}+\beta)\right]}\,.
\end{multline}
The first integral in the second line of \eqref{rewrite} converges
and therefore is a constant up to terms $O(\varepsilon)$. Asymptotic
form of the second integral for small $\varepsilon$ can be easily
found and we obtain
\begin{multline}\label{asymptotic}
    I_L(\beta)=-\frac{\pi}{4}\int_{0}^{\infty}\frac{d
    y}{y}\left\{\frac{e^{-\beta y}}{e^y-1}-\frac{1}{y\left[1+y(\frac{1}{2}+\beta)\right]}\right\}\\
    -\frac{\pi}{4}\left(\beta+\frac{1}{2}\right)\left\{\frac{1}{\varepsilon\left(\beta+\frac{1}{2}\right)}
    +\ln\varepsilon+\ln\left(\beta+\frac{1}{2}\right)\right\}+O(\varepsilon)\,.
\end{multline}
Then we rewrite $I_{\Gamma_\varepsilon}(\beta)$ making the
substitution $y=\varepsilon e^{i\varphi}$ as
\begin{multline}\label{rewrite2}
    I_{\Gamma_\varepsilon}(\beta)=\frac{1}{8}\int_{-\pi}^{\pi}d\varphi~ e^{-i\varphi}\frac{1+\left(\beta+\frac{1}{2}\right)\varepsilon
    e^{i\varphi}}{\varepsilon}\left(\ln \varepsilon+i\varphi+\frac{1}{2}\varepsilon
    e^{i\varphi}\right)+O(\varepsilon)\\
    =\frac{\pi}{4}\left(\frac{1}{\varepsilon}+
    \left(\beta+\frac{1}{2}\right)\ln\varepsilon+\frac{1}{2}\right)+O(\varepsilon)\,.
\end{multline}
Summing up $I_L(\beta)$ and $I_{\Gamma_\varepsilon}(\beta)$ and
assuming $\varepsilon\to 0$, we obtain
\begin{multline}\label{final1}
    I(\beta)=-\frac{\pi}{4}\int_{0}^{\infty}\frac{d
    y}{y}\left\{\frac{e^{-\beta y}}{e^y-1}-\frac{1}{y\left[1+y(\frac{1}{2}+\beta)\right]}\right\}\\
    +\frac{\pi}{8}-\frac{\pi}{4}\left(\beta+\frac{1}{2}\right)\ln\left(\beta+\frac{1}{2}\right)\,.
\end{multline}

Thus, we reduced the task of calculation of the initial integral to
a much easier task of calculation of the integral
\begin{equation}\label{}
    J(\beta)=\int_{0}^{\infty}\frac{d
    y}{y}\left\{\frac{e^{-\beta
    y}}{e^y-1}-\frac{1}{y\left[1+y(\frac{1}{2}+\beta)\right]}\right\}\,,
\end{equation}
which we seek as the limit
\begin{equation}\label{}
    J(\beta)=\lim\int_{0}^{\infty}\frac{d
    y}{y^{1-\alpha}}\left\{\frac{e^{-\beta
    y}}{e^y-1}-\frac{1}{y\left[1+y(\frac{1}{2}+\beta)\right]}\right\},
    \qquad \alpha\to 0\,.
\end{equation}
Since \cite{whittaker}
\begin{equation}\label{hurwitzdef}
    \int_{0}^{\infty}\frac{d
    y}{y^{1-\alpha}}\frac{e^{-\beta
    y}}{e^y-1}=\Gamma(\alpha)\zeta(\alpha,\beta+1)
\end{equation}
and \cite{whittaker}
\begin{equation}\label{}
    \int_{0}^{\infty}\frac{d
    y}{y^{2-\alpha}}\frac{1}{1+y(\frac{1}{2}+\beta)}=-\left(\beta+\frac{1}{2}\right)^{1-\alpha}\frac{\pi}{\sin\pi\alpha}\,,
\end{equation}
we obtain
\begin{equation}\label{}
    J(\beta)=\lim\left\{\Gamma(\alpha)\zeta(\alpha,\beta+1)+\left(\beta+\frac{1}{2}\right)^{1-\alpha}\frac{\pi}{\sin\pi\alpha}\right\},\qquad
    \alpha\to 0\,.
\end{equation}
Using \cite{whittaker}
\begin{equation}\label{}
    \zeta(0,\beta+1)=-\frac{1}{2}-\beta, \qquad
    {\frac{d}{d\alpha}\zeta(\alpha,\beta+1)\vline}_{~\alpha=0}=\ln\Gamma(\beta+1)-\frac{1}{2}\ln
    2\pi\,,
\end{equation}
we come to
\begin{equation}\label{j}
    J(\beta)=\gamma\left(\beta+\frac{1}{2}\right)+\ln\Gamma(\beta+1)-\frac{1}{2}\ln
    2\pi-\left(\beta+\frac{1}{2}\right)\ln\left(\beta+\frac{1}{2}\right)\,.
\end{equation}
Substitution of \eqref{j} into \eqref{final1} yields the desired
result \eqref{generalization}.

\paragraph{Corollary: proof of the formula \eqref{hurwitz}.} After
substitution $y=\ln\left(1+e^{-2ix}\right)$ in \eqref{hurwitz}, we
obtain
\begin{equation}\label{hurwitz2}
    I(\alpha,\beta)=\frac{1}{2i}\int_{-\infty}^{(0+)}\frac{e^{y(\beta+1)}}{e^y-1}y^\alpha
    dy\,.
\end{equation}
The path is the same as before. One can recognize in
\eqref{hurwitz2} the integral representation of the Hurwitz zeta
function $\zeta(\alpha+1,\beta+1)$ \cite{whittaker} (up to the
factor $-{\pi}/{\Gamma(-\alpha)}$).

\paragraph{Proof of formulae \eqref{intlog} and \eqref{intlog3}.} First, we
note that \eqref{intlog} can be presented as
\begin{equation}\label{intlog2}
    f(a)=\int_{-\frac{\pi}{2}}^{\frac{\pi}{2}}\ln\left[-a+\ln(1+e^{-2ix})\right]dx
\end{equation}
and after making the substitution $y=\ln\left(1+e^{-2ix}\right)$ as
\begin{equation}\label{equation}
    \frac{1}{2i}\int_C\ln(-a+y)\frac{e^y}{e^y-1}dy\,.
\end{equation}
The path $C$ is different for two cases i) $a<\ln 2$ and ii) $a>\ln
2$. The difference comes from the fact that the point $y=\ln 2$
should belong to the path $C$. However, there is also a cut along
the line $(-\infty, a]$. $C$ can be arbitrarily deformed such that
the pole $y=0$ is not intercepted in the course of deformation. From
this requirements, we find that in the case i) $C$ extends from
$-\infty$, circumvents the the point $y=\max\{a,0\}$, and extends
back to $-\infty$; in the case ii) it is composed of two branches:
first extends from $-\infty$ and terminates at $y=\ln 2$ at the
lower side of the cut, and the second emerges at $y=\ln 2$ at the
upper side of the cut and extends to $-\infty$. The logarithm is
$\ln(-a+y)=\ln|a-y|\mp \pi i$ at the lower and upper sides of the
cut, respectively. We also assume that $C$ does not contain the
points $y=\pm 2\pi i n$, ~$n=1,2,3,...$

Now we can write (in the case $a>0$)
\begin{multline}\label{}
    f(a)=\pi\int_\varepsilon^\infty\frac{dy}{e^y-1}-\pi\int_\varepsilon^b\frac{e^ydy}{e^y-1}\\
    +\frac{1}{2i}\int_{-\pi}^0i~d\varphi(\ln|a|+\pi
    i)+\frac{1}{2i}\int_0^{\pi}i~d\varphi(\ln|a|-\pi
    i)+O(\varepsilon)\,.
\end{multline}
Terms containing $\varepsilon$ emerge from the integration around
the pole $y=0$. Calculating the elementary integrals and proceeding
to the limit $\varepsilon\to 0$, we obtain \eqref{intlog}. If $a<0$
is the case, calculations are analogous to the presented above and
lead to the same result, therefore we omit them.

To proof \eqref{intlog3}, we present it as a sum
\begin{equation}\label{}
    g(a)=g_1+g_2\,,
\end{equation}
\begin{equation}\label{g1}
    g_1=\frac{1}{2}\int_{-\frac{\pi}{2}}^{\frac{\pi}{2}}\ln\left[-a+\ln(1+e^{-2ix})\right]e^{-2ikx}dx\,,
\end{equation}
\begin{equation}\label{g2}
    g_2=\frac{1}{2}\int_{-\frac{\pi}{2}}^{\frac{\pi}{2}}\ln\left[-a+\ln(1+e^{-2ix})\right]e^{2ikx}dx\,.
\end{equation}
\eqref{g1} and \eqref{g2} can be evaluated in the same manner as
\eqref{intlog2}. The only complication is the additional factors
$e^y-1$ and $\frac{1}{e^y-1}$ that should be embedded in
\eqref{equation}. We present the final result only:
$$
g_1=-\frac{\pi}{2}e^b
$$
$$
g_2=\frac{\pi}{2}\left(1-\frac{1}{a}+\frac{1}{e^b-1}\right)
$$

\paragraph{Proof of formulae \eqref{remarkable} and \eqref{remarkable2}.}
Although it is almost obvious how to obtain \eqref{remarkable} and
\eqref{remarkable2} from \eqref{intlog}, they are not so obvious
by themselves.  Assuming $a=\pm
n\ln\frac{1}{r}$,$~0<r<\frac{1}{2}$, $n=1,2,3,...$ in
\eqref{intlog2} and \eqref{intlog}, we obtain
\begin{equation}\label{rem1}
    \int_{-\frac{\pi}{2}}^{\frac{\pi}{2}}\ln\left[-ix-n\ln \frac{1}{r}+\ln(2\cos x)\right]dx=\pi\ln
    \left[n\ln \frac{1}{r}\right]\,,
\end{equation}
\begin{equation}\label{rem2}
    \int_{-\frac{\pi}{2}}^{\frac{\pi}{2}}\ln\left[-ix+n\ln \frac{1}{r}+\ln(2\cos
    x)\right]dx=\pi\ln\frac{n\ln \frac{1}{r}}{1-r^n}
\end{equation}
and subtracting \eqref{rem2} from \eqref{rem1}
\begin{equation}\label{n}
    \int_{-\frac{\pi}{2}}^{\frac{\pi}{2}}\ln\frac{-ix-n\ln \frac{1}{r}+\ln(2\cos
    x)}{-ix+n\ln \frac{1}{r}+\ln(2\cos
    x)}~dx=\pi\ln(1-r^n)\,.
\end{equation}
Using the following formula for the gamma function\cite{whittaker}
\begin{equation}\label{gaussproduct}
    \Gamma(x)=\lim\frac{n!n^x}{x(x+1)...(x+n)}, \qquad n\to \infty
\end{equation}
and the well known fact that\footnote{Consequence of
Lobachevskii's integral. The proof is elementary:

$\int_{0}^{\frac{\pi}{2}}\ln \cos x
dx=\int_{0}^{\frac{\pi}{2}}\ln\sin x
dx=\frac{1}{2}\int_{0}^{{\pi}}\ln\sin x dx=\frac{\pi}{2}\ln
2+2\int_{0}^{\frac{\pi}{2}}\ln\cos x dx$}
\begin{equation}\label{lobachevsky}
    \int_0^{\frac{\pi}{2}}\ln(2\cos x)~dx=0\,,
\end{equation}
one finally comes to \eqref{remarkable}.

\eqref{remarkable2} is obtained from (see \eqref{intlog2} and
\eqref{intlog})
$$
\int_{-\frac{\pi}{2}}^{\frac{\pi}{2}}\ln\left[-ix-n-c+\ln(2\cos
    x)\right]~dx=\pi\ln(n+c)\,,\quad c>\ln 2\,,\quad
    n=0,1,2,...
$$
with the aid of identities \eqref{gaussproduct} and
\eqref{lobachevsky}.

\section{Some other results}
We end this article by listing some integrals without proofs, which
are either the direct consequence of the integrals in the main text
or new integrals, which can be evaluated exactly in the same manner.
\begin{equation}\label{}
    \int_{0}^{\frac{\pi}{2}}\ln[{x^2+\ln^2\cos x}]~d
    x=\frac{\pi}{2}\ln\ln2
\end{equation}
\begin{equation}\label{}
    \int_{0}^{\frac{\pi}{2}}\ln[{x^2+\ln^2\cos x}]\cos 2x~d
    x=-\frac{\pi}{\ln2}
\end{equation}
\begin{equation}\label{}
    \int_{0}^{\frac{\pi}{2}}\ln[{x^2+\ln^2(2\cos x)}]\cos 2x~d
    x=-\frac{\pi}{4}
\end{equation}
\begin{equation}\label{}
    \int_{0}^{\frac{\pi}{2}}\frac{\ln\cos x}{x^2+\ln^2\cos x}~d
    x=\frac{\pi}{2}\left(1-\frac{1}{\ln 2}\right)
\end{equation}
\begin{equation}\label{}
    \int_{0}^{\frac{\pi}{2}}\frac{x\sin 2x}{x^2+\ln^2\cos x}~d
    x=\frac{\pi}{4\ln^22}
\end{equation}
\begin{equation}\label{}
    \int_{0}^{\frac{\pi}{2}}\frac{x\sin 2x}{x^2+\ln^2(2\cos x)}~d
    x=\frac{13\pi}{48}
\end{equation}
\begin{equation}\label{}
    \int_{-\frac{\pi}{2}}^{\frac{\pi}{2}}\frac{\left(1+e^{-2ix}\right)^\beta}{\ln\left(1+e^{-2ix}\right)-a}~d x
    =-\frac{\pi}{a}+\pi\frac{e^{(\beta+1)a}}{e^a-1}H(\ln 2-a)
\end{equation}
\begin{equation}\label{}
    \int_{0}^{\frac{\pi}{2}}\frac{x\sin 2x}{x^2+\ln^2(2e^{-a}\cos x)}~d
    x=\frac{\pi}{4a^2}+\frac{\pi
    e^{a}}{4}\left(1-\frac{1}{(e^a-1)^2}\right)H(\ln2-a)
\end{equation}



\end{document}